# Computing Feedback Laws for Linear Systems with a Parallel Pieri Homotopy

Jan Verschelde
Dept of Math, Stat, and CS
University of Illinois at Chicago
851 South Morgan (M/C 249)
Chicago, IL 60607-7045, USA.
Email: jan@math.uic.edu
URL: http://www.math.uic.edu/~jan

Yusong Wang
Dept of Math, Stat, and CS
University of Illinois at Chicago
851 South Morgan (M/C 249)
Chicago, IL 60607-7045, USA.
Email: ywang25@uic.edu
URL: http://www.math.uic.edu/~ywang25

*Abstract*— **Homotopy methods to solve polynomial systems are well suited for parallel computing because the solution paths defined by the homotopy can be tracked independently. Both the static and dynamic load balancing models are implemented in C with MPI, adapting PHCpack written in Ada using gcc, and tested on academic benchmarks and mechanical applications. We studied the parallelization of Pieri homotopies to compute all feedback laws to control linear systems. To distribute the workload, we mapped the poset onto a tree. As the dimensions of the Pieri homotopies grow incrementally from the root to the leaves in the tree, we found the Pieri homotopies well suited for parallel computing.**

## I. INTRODUCTION

Polynomial systems occur in a wide variety of application domains, such as mechanical design, signal processing, and in what is most relevant for this paper: the control of linear systems. Typically, the number of solutions of a polynomial system grows exponentially with its dimension. For example, the polynomial system whose solutions are feedback laws to control a machine with $m$ inputs and $p$ outputs has in general as many solutions as $d_{m,p} = \frac{1! 2! 3! \cdots (p-2)!(p-1)! \cdot (mp)!}{m!(m+1)!(m+2)! \cdots (m+p-1)!}$. So the need for parallel computation is very real.

Homotopy continuation methods are reliable and powerful methods to compute numerical approximations to all isolated complex solutions. These methods operate in two stages. To solve a system $f(\mathbf{x}) = \mathbf{0}$, we first construct $g(\mathbf{x}) = \mathbf{0}$ whose solutions are known. This system $g(\mathbf{x}) = \mathbf{0}$ then becomes the start system in the homotopy

$$h(\mathbf{x}, t) = \gamma g(\mathbf{x})(1-t) + t f(\mathbf{x}) = \mathbf{0}, \quad \gamma \in \mathbb{C}. \quad (1)$$

For almost all choices of the complex constant $\gamma$, all solutions paths $\mathbf{x}(t)$ are regular and bounded for $t \in [0, 1)$. In the second stage, continuation methods are applied to track the paths starting at the known solutions of $g(\mathbf{x}) = \mathbf{0}$ to the desired solutions of $f(\mathbf{x}) = \mathbf{0}$, as $t$ goes from 0 to 1. See [11] for a survey on recent methods to construct efficient homotopies.

The publicly available software PHCpack [20] implemented the homotopy methods in a sequential version.

The homotopy algorithm is well suited for parallel computing, since the paths can be tracked independently from each other. The efficiency of the algorithms for solving systems of nonlinear equations using probability-one homotopy methods in parallel is discussed in [1], [3], [7]. More recently, in [6] and [19] the authors report on a parallel implementation of polyhedral homotopy methods, which exploit the sparse structure of polynomial systems.

In this paper we first outline the extension of PHCpack with a parallel path tracker, before describing the parallelization of Pieri homotopies to solve the pole placement problem in the control of linear systems. This paper is a sequel to [21] and [22].

April 15, 2004

## II. A PARALLEL PATH TRACKER IN PHCPACK

In this section we describe two load balancing schemes and report computational experiences on tracking paths defined by the homotopy (1), for two large polynomial systems.

### A. Static and Dynamic Workload Balance

For best performance, the workload should be distributed evenly among the processors. In the static workload distribution, the paths are distributed evenly to the processors once at the start. While this leads to a minimal communication overhead, the workload for each processor may have a large variance, as paths diverging to infinity require more time. The dynamic workload assignment with a master/slave paradigm is usually better. Each of the slave processors will be given one job at the beginning. After a slave finishes its job, it sends the result to the master, which sends then a new job to the slave. While this requires more communication overhead than the static workload assignment model, we can improve it by overlapping the communication and computation with the non-blocking sending and receiving in the MPI library.

## B. Experimental Results and Discussion

Our parallel code was developed on a rocketcalc atlas cluster with four 2.4 GHz processors under Linux. To examine the speedup and the load balancing issues better on larger problems, we ran the code on the Platinum cluster at NCSA.

*1) An Academic Benchmark: cyclic 10-roots:* The cyclic $n$-roots problem is widely used as a benchmark for publicly available software ([5], [6], [20]). Computing all cyclic $n$-roots is hard because the number of paths is often too large to be traced by a single computer [4]. For $n = 10$, we need to trace 35,940 paths. With a given start system, it takes 8 hours with the sequential version of path tracker on a 1GHz computer. Our parallel path tracker traces all 35,940 paths within 5 minutes on 128 1GHz CPUs.

TABLE I

SPEEDUPS OF THE STATIC AND DYNAMIC LOAD BALANCING FOR THE CYCLIC 10-ROOTS PROBLEM ON THE PLATINUM CLUSTER AT NCSA. TIME UNITS ARE USER CPU MINUTES.

| #CPUs | Static | | Dynamic | | Improvement dynamic/static |
|---|---|---|---|---|---|
|  | time | speedup | time | speedup |  |
| 1 | 480.0 | 1.0 | 480.0 | 1.0 | – |
| 8 | 75.5 | 6.4 | 66.6 | 7.2 | 11.75% |
| 16 | 36.4 | 13.2 | 31.7 | 15.2 | 12.96% |
| 32 | 19.0 | 25.3 | 15.7 | 30.7 | 17.56% |
| 64 | 10.2 | 46.9 | 7.9 | 60.5 | 22.48% |
| 128 | 6.6 | 73.3 | 4.3 | 112.9 | 35.11% |

From Table I, we see that the dynamic workload balancing improves the total time of the static approach by 10% to 35%. For this problem, the variance of the time needed to trace the paths can be large (one thousand paths diverge). The improvement of using dynamic load balancing is more obvious with more processors since the variance becomes larger for fewer jobs on each processor in the static workload assignment. Fig. 1 shows that the speedup is almost optimal for the dynamic model when the number of processors is less than 32. For any number of processors, dynamic load balancing wins.

*2) An Application from Mechanism Design:* This example comes from the geometric design of the five degree-of-freedom robot formed by links connected by revolute, prismatic and spherical joints to form an RPS serial chain [16], [17], [18]. To design this robot one must solve ten polynomial equations in ten unknowns. The homotopy we used (using a linear-product start system as in [17]) led to 9,216 solution paths. As reported in [16], the sequential version of PHCpack takes about 24 hours on a 2.4GHz Pentium IV machine. On 128 1GHz CPUs of the Platinum cluster at NCSA, all paths were traced within 22 minutes. As the time for one 1GHz was unavailable, we assumed an initial optimal speedup in the dynamic case, extrapolating to 3111.2 CPU minutes sequential time, obtained as $8 \times 388.9$, see Table II. While assuming an initial optimal speedup with 8 CPUs is unrealistic, when doubling the number of processors we may finish more than twice as fast when the distribution of the workload is more evenly spread among the processors, see Fig. 2.

TABLE II

COMPARISON OF THE STATIC AND DYNAMIC WORKLOAD BALANCE WITH THE RPS PROBLEM. TIME UNITS ARE USER CPU MINUTES.

| #CPUs | Static | | Dynamic | | Improvement dynamic/static |
|---|---|---|---|---|---|
|  | time | speedup* | time | speedup* |  |
| 8 | 417.5 | 7.5 | 388.9 | 8.0 | 6.84% |
| 16 | 195.1 | 15.9 | 183.7 | 16.9 | 5.84% |
| 32 | 94.7 | 32.9 | 96.1 | 32.4 | -1.50% |
| 64 | 49.8 | 62.5 | 47.5 | 65.5 | 4.65% |
| 128 | 25.1 | 124.0 | 22.0 | 141.4 | 12.43% |

In Table II, we can find the improvement of the dynamic over the static balancing model is not obvious here. Since in this example, there are more than eight thousand diverging paths, which dominate the total computation time and each of the diverging path spend almost the same time. So there is no large variance in the workload among the processors in the static model. Moreover, the overhead of the communications decreases the efficiency of the dynamic load balancing model.

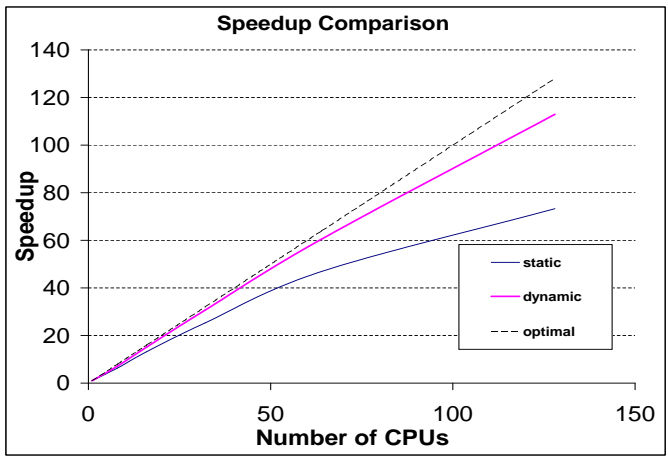

Fig. 1. The speedup comparison of the static and dynamic load balancing for the cyclic 10-roots problem on the Platinum Cluster at NCSA.

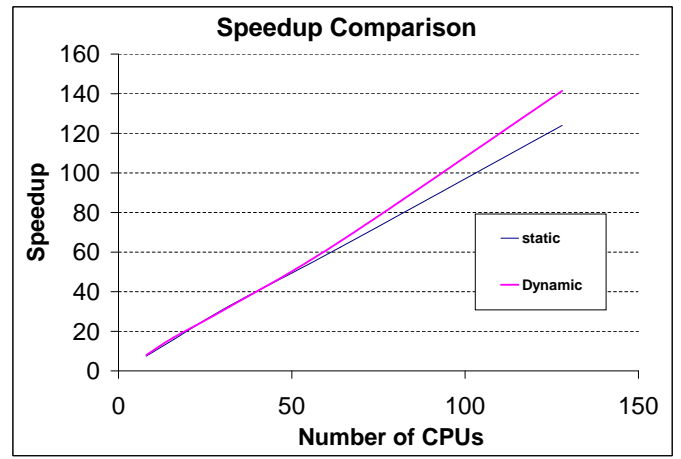

Fig. 2. The speedup comparison of the static and dynamic load balancing for the mechanical application on the Platinum Cluster at NCSA.

We took this application to illustrate what happens when – an unfortunately still too often occurring case – many solution paths diverge to infinity. For this particular system, the mixed volume gives the exact root count of 1024, and thus the polyhedral homotopy (implemented in [5], [6], and [20]) is optimal. The black-box solver of PHCpack gives the complete solution list in 24.6 minutes CPU time on a 2.4Ghz Linux machine. The system and its solutions are available online at http://www.math.uic.edu/~jan/demo.html.

## III. PARALLEL PIERI HOMOTOPY

In this section, we consider the parallel implementation of a homotopy continuation method to find a general start system $g(\mathbf{x}) = \mathbf{0}$ to be used in the homotopy (1) to solve a particular problem $f(\mathbf{x}) = \mathbf{0}$.

### A. Solving the Pole Placement Problem

As in the introduction, we consider a machine with $m$ inputs and $p$ outputs, whose evolution in time is governed by a system of linear differential equations. The control of this machine by a compensator with $q$ internal states corresponds to a problem from enumerative geometry for which the so-called Pieri homotopies were derived. The theoretical connection was first made in [2] for $q = 0$, and generalized in [13], [14], and [15]. Algorithms, defined by Pieri homotopies were developed in [8], [9] and [12]. In this paper we explain the algorithms geometrically. For the relation with inverse eigenvalue and matrix extension problems, see [10].

For a given $m$, $p$, and $q$, denote $n = mp + q(m + p)$. The problem we solve takes on input $n$ general $m$-planes $L_i$ in $\mathbb{C}^{m+p}$ and $n$ interpolation points $s_i \in \mathbb{C}$, $i = 1, 2, \ldots, n$. For this input, we want to compute all polynomial maps $X(s)$ of degree $q$ producing $p$-planes that meet those given general $m$-planes $L_i$ at the prescribed interpolation points $s_i$, i.e.: we are given $n$ intersection conditions:

$$\det(X(s_i)|L_i) = 0, \quad i = 1, 2, \ldots, n. \quad (2)$$

These intersection conditions define a polynomial system in the coefficients of the map $X : \mathbb{C} \to \mathbb{C}^{(m+p) \times p} : s \mapsto X(s)$.

$$X(s,t) = \begin{bmatrix} \star & 0 + \star s \\ \star & \star t + \star s \\ \star & \star t + \star s \\ \star & \star t + 0s \end{bmatrix} \Leftrightarrow 4 \to \begin{bmatrix} \star & 0 \\ \star & \star \\ \star & \star \\ \star & \star \\ 0 & \star \\ 0 & \star \\ 0 & \star \leftarrow 7 \\ 0 & 0 \end{bmatrix} \Leftrightarrow \begin{matrix} [4\ 7] \\ \text{Shorthand} \end{matrix}$$

Standard      Concatenated

Fig. 3. Localization pattern of solutions for $p$=2, $m$=2, $q$=1

In what follows, we show that the problem (2) is well posed: we have $n$ equations in the $n$ variables which define a general map $X(s)$.

### B. Localization Patterns and Pieri Homotopies

We represent $X(s)$ by a localization pattern in $\{0, \star\}^{(m+p) \times p}$ (i.e.: a matrix over $\mathbb{Z}_2$) in which all stars stand for the nonzero coefficients of the generator matrix. A $p$-plane fits a localization pattern if it can be represented by a matrix of generators with zero entries everywhere the localization pattern prescribes them. For example, in Fig. 3, the left picture is the canonical form of the degree one-map solution localization pattern for $p = 2$, $m = 2$, $q = 1$, where the $t$ is for homogenizing the polynomials to deal with both bottom pivots and top pivots. The middle picture is the concatenated form with $p + n$ stars, where we append the higher degree coefficients below the lower degree coefficients and the degree of freedom is $n = 8$. The right picture is a shorthand notation for the bottom pivots which record the row indices of the bottommost stars.

A valid bottom pivot localization pattern is defined as below:

1) Let $q = dp + r$ with $d, r \in \mathbb{N}$ and $r < p$. A localization pattern for $(m + p) \times p$-maps of degree $q$ has the first $p - r$ columns with dimension $(d + 1)(m + p)$ and the remaining columns have dimension $(d + 2)(m + p)$.
2) All stars within a column should be contiguous and the row indices in which the bottommost and topmost stars occur strictly increasing as a function of the column index. These indices are called the top and bottom pivots, respectively.
3) No two bottom pivots differ by $m + p$ or more.

The above definition is extracted from [9]. In this paper – and in our preliminary parallel implementation – we consider the top pivots as fixed to $[1\ 2\ \cdots\ p]$.

The special emphasis on the format of $X(s)$ is entirely justified as it leads naturally to a homotopy as follows. The bottommost pivots of $X$ give a recipe (see [9]) for a special $m$-plane $S_X$ so that $\det(X|S_X) = 0$ if and only if at least one of the entries in $X$ at the bottommost pivots is zero. The *Pieri homotopy* in (3) moves $S_X$ to $L_n$. For $q > 0$, the map $X(s)$ can meet $S_X$ only at $\infty$, so the corresponding interpolation point moves from $\infty$ to $s_n$. To represent $\infty$ properly, we homogenize the polynomials in $X(s)$ using $t$, and denote the maps as $X(s, t)$. The Pieri homotopy in (3) then moves $(s, t) = (1, 0)$ to $(s_n, 1)$. Observe the double use of $t$ in (3): as continuation parameter and variable added to homogenize the maps.

$$H(X(s,t), s, t) =$$
$$\begin{cases} \det(X(s,t)|(1-t)S_X + tL_n) = 0 \\ (s-1)(1-t) + (s-s_n)t = 0 \\ \det(X(s_i,t_i)|L_i) = 0, \quad i=1,2,\ldots,n-1 \end{cases} \quad (3)$$

for $t \in [0, 1]$.

The start solutions for the Pieri homotopy all fit in the patterns $Y(s,t)$ obtained from $X(s,t)$ by turning a bottommost star to zero. By induction on $n$, we assume that all these children $Y(s,t)$ meet already the $n-1$ general $m$-planes $L_i$ at $(s_i, t_i)$, i.e.: $\det(Y(s_i,t_i)|L_i) = 0$, for $i = 1, 2, \ldots, n-1$. To satisfy the $n$th intersection condition with $L_n$ at $(s_n, t_n)$, we trace the solution paths defined by the Pieri homotopy in (3), as t goes from 0 to 1.

Note that at $t = 1$, we have $s = s_n$, and all intersection conditions in (2) are satisfied. In the next section we describe the induction on $n$ which leads to an efficient way to count all the roots.

### C. Counting Roots by Posets and Trees

The shorthand notation of the bottom pivots is a convenient way to count the solution maps and to represent the nested sequences of homotopies needed to compute all solutions.

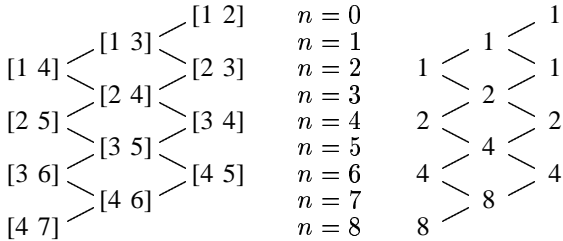

Fig. 4. Combinatorial root count for $p = 2$, $m = 2$, $q = 1$ with the poset structure. The brackets at the above are the bottom pivots. The trivial localization pattern corresponds to [1 2]. The solutions are counted at the below, starting at the top and adding up the numbers at the leaves while moving down to the root of the poset [4 7] yielding 8 solutions.

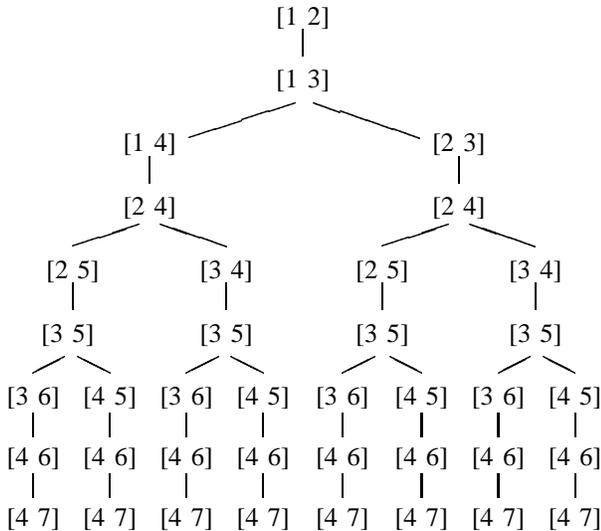

Fig. 5. Combinatorial root count for $p = 2$, $m = 2$, $q = 1$ with the Pieri tree. The brackets at the above are the bottom pivots. The trivial localization pattern corresponds to [1 2]. The solutions are counted by starting at the top and finding all the allowable paths to reach the leave [4 7], adding up the number of leaves through different paths yielding 8 solutions.

In Fig. 4, the poset structure is described to count the number of the solution planes. It starts from the trivial localization pattern which has its top pivots [1 2 $\cdots$ p], e.g.: [1 2], and its bottom pivots [$d(m+p)+m+r+1$ $\cdots$ $d(m+p)+m+p$ ($d+1)(m+p)+m+1$ $\cdots$ $(d+1)(m+p)+m+r$], e.g.: [4 7]. The bottom poset structure at the left in Fig. 4 is obtained by fixing the top pivots and decreasing a bottommost pivot, which is called a bottom child, to get all the valid localization patterns recursively. At the right of Fig. 4 we see the counting procedure: we start with one solution for $n = 0$; then, for $n > 0$, the number of maps fitting $X$ and meeting $n$ general $m$-planes equals the sum of the number of solution maps fitting the children of $X$ and meeting $n - 1$ general planes. Every link in the poset as in Fig. 4 corresponds to one instance of the Pieri homotopy.

The combinatorial root count with a poset is implemented in the sequential version of PHCpack [20]. In the parallel version of the Pieri homotopy program, we solve the problem based on a tree, called Pieri trees in [8]. The tree corresponding to the poset in Fig. 4 is shown in Fig. 5.

To see the virtue of Pieri trees for parallel computers, we need to recall the induction on $n$ in the derivation of the Pieri homotopies. In the Pieri tree, each edge represents a job, i.e.: the tracking of one solution path. Two jobs (represented by two edges in the Pieri tree) become completely independent from each other once the solution at their common ancestor node has been computed. Compared to posets, the organization of the path tracking along Pieri trees makes them more suitable for parallel computers, since the workload for each of the processors can be balanced well.

Also the memory management becomes simpler with Pieri trees. Every node in the tree is only needed in the computation of the path ending at the node, or in the paths originating at the node, in total in no more than $p+1$ jobs. So in general, the memory occupied by a node in the Pieri tree can be released rather quickly after a job has finished. In the poset however, the nodes carry the information of many more paths and need to remain active even if only one job is still not completed. Especially for larger problems, as the number of roots grows exponentially, the number of internal nodes may also increase dramatically, exhausting all the memory rather quickly.

As pointed out earlier, the Pieri homotopies used in this paper keep their top pivots fixed, while the full blown sequential version of the program allows to increase top and decrease bottom pivots simultaneously, hereby satisfying two new intersection conditions with one Pieri homotopy. As this scheme needs in general fewer solution paths than when keeping top pivots fixed, we plan to extend our preliminary parallel implementation in the near future to the most general Pieri homotopies.

### D. Parallel Pieri Homotopy Algorithm

Fig. 6 illustrates the procedure of the parallel Pieri computation. Here we apply the dynamic workload assignment with a master/slave paradigm, which is proved to be better in the previous section. At the beginning, the master generates (at most $p$) jobs by increasing the bottom pivots and put them in the queue. Then the master distributes the available jobs to the slaves, which will finish the computation task. When one of the slave finishes its job, it returns the result to the master. The master generates (at most $p$) new jobs according to the returned information, which includes the pivot information of the node and the target solution of the previous homotopy. At first, the jobs are distributed sequentially according to the rank of the slave. After all of the slave processors are activated, the dynamic workload balance paradigm based on the first-come-first-serve strategy is implemented to compute the remaining jobs.

TABLE III
NUMBER OF PATHS AND USER CPU TIMES FOR $m=2$, $p=2$, AND $q=1$.

| $n$ | #paths | user CPU time |
|---|---|---|
| 1 | 1 | 0ms |
| 2 | 2 | 0ms |
| 3 | 3 | 10ms |
| 4 | 5 | 30ms |
| 5 | 8 | 80ms |
| 6 | 13 | 370ms |
| 7 | 21 | 1s 290ms |
| 8 | 34 | 3s 830ms |
| 9 | 55 | 8s 190ms |
| 10 | 55 | 7s 840ms |
| 11 | 55 | 16s 570ms |
| Total | 252 | 38s 350ms |

Since the workload is dynamically distributed, the slave does not know the number of jobs needs to be done in advance. So we need figure out a way to terminate the computation subroutine properly. The intuitive idea is when all of the slaves returned a leave, which can not generate any new jobs, we are done. It could be the case, when some of the slaves return leaves and find no job in the queue, they won't work any more, while the other slaves are still working on the internal nodes of the virtual tree. This situation will cause an unbalanced workload distribution, therefore the efficiency will be lower. To avoid this case, we maintain another queue to record which slave has returned a leaf, and activate it again when there are more jobs available to compute. After all the slaves have returned leaves, the master will send a message to each slave to terminate their busy waiting loop.

As the path tracking jobs are subject to a tree hierarchy, every job has to wait till the first path starting at the root node has terminated. Every job in the tree has to wait till the job providing its start solution has finished. So at the start of the of program, only very few processors are active, while most other processors are idle, waiting for their start solutions. Fortunately, the jobs closest to the root are the smallest, as exemplified in Table III. Typically, almost half of the time is spent at the last level, towards the leaves of the Pieri tree.

## IV. APPLICATIONS

Table IV shows the experimental results for different values of $m$, $p$, and $q$. As $m$, $p$, and $q$ increase, the number of solutions grows exponentially, for example: 135,660 for $p=4$, $m=3$, $q=1$. The dimension $n$ of the problem grows too, but fortunately as a polynomial, i.e.: $n = mp + q(m+p)$. Nevertheless, the problem of computing **all** solutions quickly becomes intractable.

Our 2.4GHz PC under Linux can only solve some low dimensional problems in hours. On the Platinum cluster at NCSA, we improved the time from hours to minutes for some lower dimensional problems and solved some higher dimensional problems which are not tractable for our PC. The table is drawn in an upper triangular format to show the limit of the problem which a PC can solve.

## ACKNOWLEDGMENTS

We thank the National Center for Supercomputing Applications (NCSA) for the use of the Platinum IA32 Cluster. This material is based upon work supported by the National Science Foundation under Grant No. 0105739 and Grant No. 0134611.

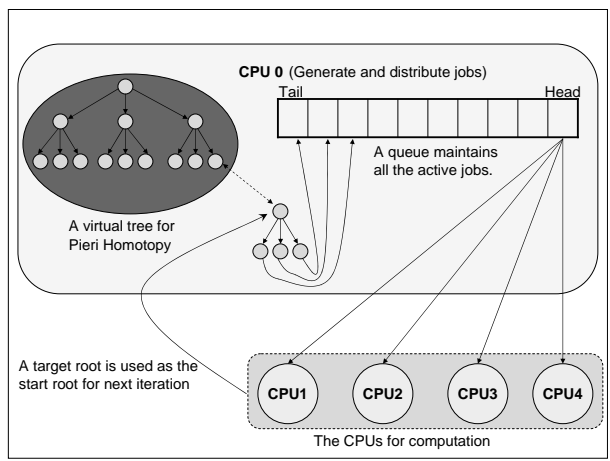

Fig. 6. Parallel Pieri homotopy with a virtual tree structure

TABLE IV

SOLVING PIERI HOMOTOPY PROBLEM ON A 2.4GHZ PC AND 64 1GHZ CPUS OF PLATINUM CLUSTER AT NCSA. TIME UNITS ARE USER CPU SECONDS.

| $p$ | $m$ | $q=0$ | | | $q=1$ | | | $q=2$ | | | $q=3$ | | |
|---|---|---|---|---|---|---|---|---|---|---|---|---|---|
| | | #Solutions | time(s) PC | Cluster | #Solutions | time(s) PC | Cluster | #Solutions | time(s) PC | Cluster | #Solutions | time(s) PC | Cluster |
| 2 | 2 | 2 | 0.2 | – | 8 | 0.9 | – | 32 | 18.4 | – | 128 | 218.3 | 19.1 |
| 3 | 2 | 5 | 0.2 | – | 55 | 38.4 | – | 610 | 2331.7 | 137.2 | 6765 | N/A | 4749.0 |
| 3 | 3 | 42 | 8.8 | – | 2730 | 7663.8 | 327.7 | 17462 | N/A | – | | | |
| 4 | 3 | 462 | 638.7 | 52.4 | 135660 | N/A | – | | | | | | |
| 4 | 4 | 24024 | N/A | 1891.2* | | | | | | | | | |

*done on 256 CPUs